# Symmetries in Images on Ancient Seals

**Amelia Carolina Sparavigna**

Dipartimento di Fisica, Politecnico di Torino, Corso Duca degli Abruzzi 24, 10129 Torino, Italy
E-mail: amelia.sparavigna@polito.it

**Abstract:** We discuss the presence of symmetries in images engraved on ancient seals, in particular on stamp seals. Used to stamp decorations, to secure the containers from tampering and for owner's identification, we can find seals that can be dated from Neolithic times. Earliest seals were engraved with lines, dots and spirals. Nevertheless, these very ancient stamp seals, in the small circular or ovoid space of their bases, possess bilateral and rotational symmetries. The shape of the base seems to determine the symmetries of images engraved on it. We will also discuss what could be the meaning of antisymmetry and broken symmetry for images on seals.

**Keywords:** Symmetry, Seals, Antisymmetry.

## 1. Introduction

We can find symmetries in human artifacts across time and cultures. As widely discussed, this is considered as due to a biological significance of symmetry [1,2]. The simplest symmetry, the bilateral one, seems to be deeply connected with the human perception of health and then beauty of living beings. Humans are very sensitive to the presence of symmetry in images [3] and many speculations about the mechanisms underlying the vision of symmetry have been proposed as proposed [4,5].
This human sensitivity to symmetries not only gives rise to an aesthetic sense of pleasure, when good proportioned and balanced forms are observed [6], but also stimulates the use of symmetries in science [1]. In physics symmetry means the invariance of the physical laws under specific transformations. Conservation of energy, momentum and angular momentum can be viewed as a consequence of symmetry for continuous translations in time and space and for rotation, according to the Noether's Theorem. In quantum mechanics, wave functions can be either symmetric or antisymmetric, depending on particles' spin. Particles with antisymmetric wave functions are called fermions and obey the Pauli exclusion principle. Particles with integer spin have a symmetric wave function and are called bosons. In solid state physics, discrete symmetries govern oscillations of the atomic structure and mechanisms of charge and energy transport.

Studies on the symmetries in pottery decorations have been already performed. Important to identify the period of a certain culture, symmetries have a practical use in the restoration of damaged materials. In this paper, we will discuss the presence of symmetry in other ancient artefacts, the stamp seals. Considered as a form of art of non-primary importance, engraved images can be not only beautiful but also rather important for archaeological studies, because, as pottery, they are durable artefacts.

Ancient stamp seals and cylinder seals were made of a hard material and used to press an engraved figure into soft clay, for securing purposes [8], or, as for early Neolithic seals, to stamp decorations on leather and textile. It is possible that they were used to stamp non-permanent tattoo, as guessed for seals found at Çatalhöyük (7500-5700 BCE) [9]. This use is guessed also for pre-Columbian seals [10].

Cylinders firstly appeared in Mesopotamia in the late 4th millennium BCE, then widespread in Syria and Egypt and in the Aegean area. Stamp seals preceded cylinders. Early stamp seals were used in Persia, northern Syria, and Anatolia [8]. In Egypt, the scarab seals largely replaced cylinder seals early in the 2nd millennium BCE and continued as the main type, till they were replaced by the signet ring in Roman period. For Egyptians, the scarab seal was not only an impression seal: it was also an amulet with images and symbols engraved to protect the owner [11-14]. Egyptian scarabs have inscriptions, human and animals figures. Many scarabs of the Middle Kingdom Period have the base decorated with coils and entanglements of cords: these seals display bilateral and two-fold rotational symmetries (see Fig.1).

Seals collections are seldom published in public domain, then it is impossible to be exhaustive in describing symmetries displayed by these objects. In this paper, which is just an introduction to more detailed researches, let us try to follow a time-line to see if any evolution in decoration symmetries was present and what are the constrains imposed by the shape of the seal.

**Figure 1.** Seals of the Middle Kingdom with spiral coils and crosses, at the Egyptian Museum of Torino.

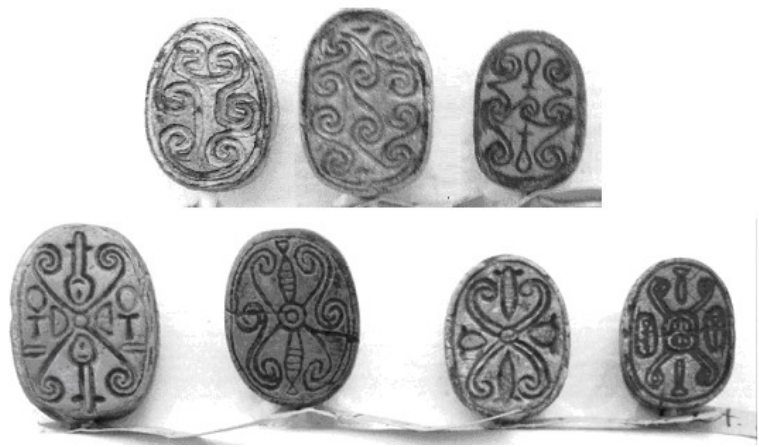

## 2. Symmetries and seals

Images on stamp seals can display symmetries adequate with the shape of their bases. Mirror or bilateral symmetry occurs when the two halves of a whole are each other's mirror images. Rotational symmetries occur with respect to rotations in the space: rotations are about an axis perpendicular to the plane of the image. We find a two-fold symmetry if we rotate the image of 180 degrees, a four-fold one when the rotation is of 90 degrees. Generally, we have *n*-fold rotation when angle is $2\pi/n$ .Because it is quite natural to follow the image in its rotation in a round space, the rotational symmetry produces a feeling of motion and evolution.

We can ask ourselves how old is the presence of rotational symmetries in images. Besides pottery, the other durable objects that can help us in investigating this question are the stamp seals. Fig.2 shows some stamp seals from Çatalhöyük (7500-5700 BCE). This was a very large Neolithic and Chalcolithic settlement in southern Anatolia. Moreover, this is the best preserved Neolithic site found to date. Images have geometric patterns with bilateral and rotational symmetries. Fig.2 shows also an interesting seal with an image which could be a very primitive symbol for interconnection of complementary opposites. A four-fold symmetry is displayed by a round stamp seal [9].

**Figure 2.** Seals from Çatalhöyük [9] have geometric images with symmetries: (a) seems to have a two-fold rotation symmetry, as seals (b),(e) and (f). (c) has a bilateral symmetry. (d) tends to a four-fold symmetry. Stamp seals with spiral pattern have also been found at Çatalhöyük (g).

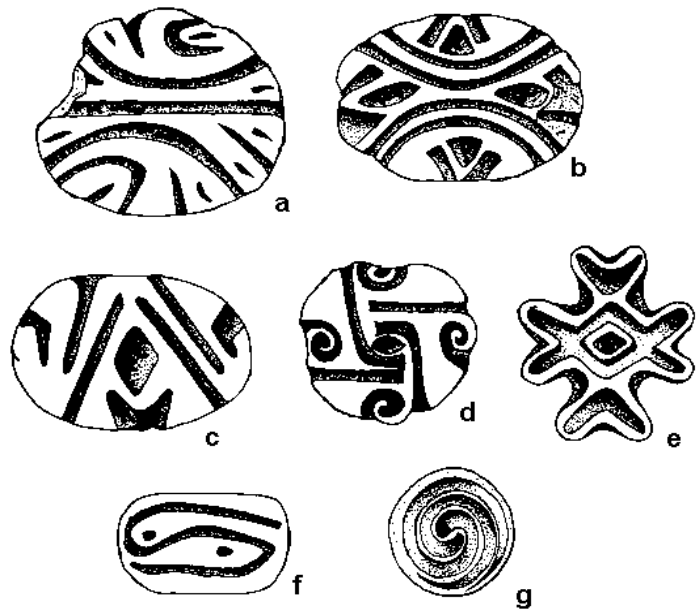

The impressing of carved stones into clay to seal containers had a long tradition in Near East, in particular in Mesopotamia, with the earliest evidence found in Syria dating to the seventh millennium BCE. During the Ubaid period, the variety of designs carved on seals expanded from the simple geometric forms to include animals with humans, snakes and birds. Seals decorated with four-legged horned animals can be easily found.

In a stamp seal from the region North Syria - Iraq [15], dated 5th-4th millennium BCE, we see a standing male figure between two horned quadrupeds back to back and head to end (see Fig.3a). The overall structure of the image is built to respect the two-fold rotational symmetry of the animals.

**Figure 3.** Images on a stamp seal from the region North Syria- Iraq (image adapted from Ref.[15]), dated 5th-4th millennium BCE (on the left) and on a stamp from Susa (image adapted from Ref.[16]), at Louvre Museum (on the right).

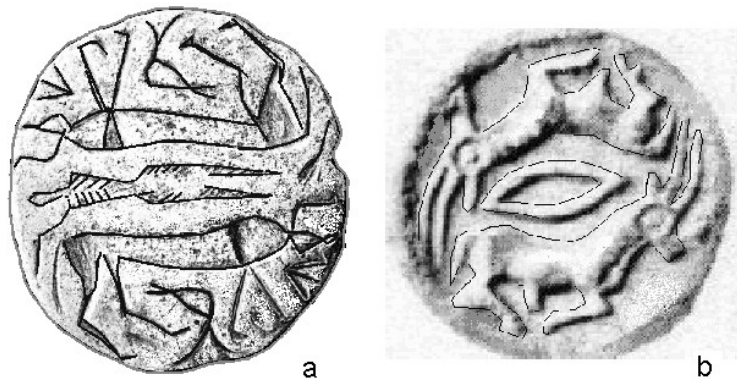

At its right part, Figure 3 (b) shows a stamp seal from Susa, in Iran. Susa is one of the oldest known settlements of the world, possibly founded about 4200 BC, although the first traces of an inhabited village have been dated to ca. 7000 BC. The seal depicts two goat-antelopes head to tail, outside an oval [16]. In this seal, the idea of motion is strongly enhanced, the two antelopes seem to run on the rim of the seal. At the same Ref.16, other seal images are shown: unfortunately, the quality of images is very low, but nevertheless it is possible to argue that they have a bilateral symmetry, in particular the seals from Dilmun culture. Dilmun, the land mentioned by Mesopotamian as a trade partner and as a trading post of the route to the Indus Valley, might be associated with the islands of Bahrain. Figure 4 shows a Dilmun seal, the left part of which has a mirror symmetry, with stylistic influences from the Indus Valley culture.

Stamp seals of the ancient Indus cities are quite important, because they bring inscriptions in a writing system which is not yet deciphered. There is another fact which is quite interesting to note. These seals are usually engraved with single human or animal figures. It is then very difficult to find images displaying symmetries (one is shown in Fig.4), and it could be rather interesting to find a reason from this lack of symmetry. Nevertheless there are small Indus Valley seals with four-fold symmetries, showing crosses and swastikas, suitable with the square shape of the seals.

As told in the introduction, symmetry is a common phenomenon in human culture. This is easy to understand for bilateral symmetry which seems to be quite naturally associated with the fact that animals and human beings possess this symmetry. But other symmetries, which are rare in nature, for instance the two-fold rotational symmetry, are quite common in very old artefacts. A possible reason for the appearance of rotational symmetries in decorations was proposed by P. Gendes. According to him [17,18], these symmetries can naturally arise in solving problems involved in weaving activities. It is then not surprising that the earliest seals have geometric designs with rotational symmetries.

**Figure 4.** Stamp seal from Dilmun (a) and Indus Valley seals (b) and (c).

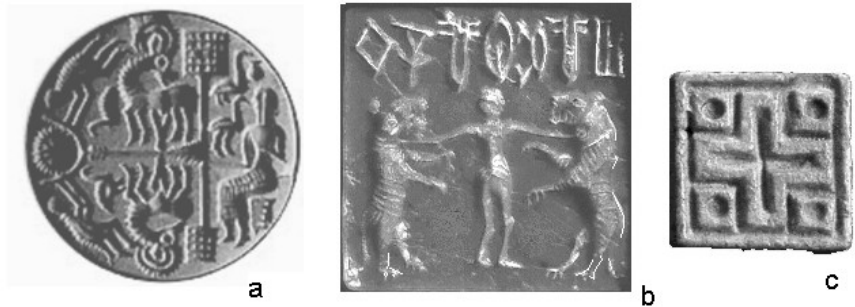

## 3. Antisymmetries and broken symmetries.

An antisymmetry in images means that we observe the mirror or the rotated image in the opposite colour. When the image has a rotational antisymmetry, we have an enhancement of the feeling of motion and evolution, due to the contrast of interpenetrating colours. Let us consider for instance a well-know symbol, the Yin-Yang symbol of Taoism. As in many religious symbols, the symmetry is used to convey an intuitive meaning. And in fact, we see that the ancient Taijitu image of Taoism has a fascinating use of symmetry for rotation about the central point (a two-fold one), combined with black-and-white inversion of colours (see Fig.5). The image actually intends to be a representation of the complementary need for male and female concepts. The rotational antisymmetry catches our attention with a feeling of development in the frame of complementary actions, not of mere recurrence.

**Figure 5.** A two-fold rotation symmetry on the left gives an idea of recurrence. The anti-rotational symmetry is shown by the Yin-Yang symbol (on the right).

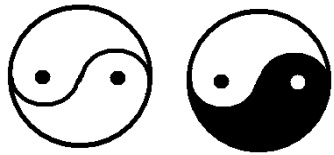

The following image shows how the use of bilateral antisymmetry can create interesting shadowing and three-dimensional effects. A use of antisymmetry can create beautiful effects in the body art, as the picture at the right side of Fig.6 is showing.

**Figure 6.** Bilateral antisymmetry can create a shadowing effect. Opposite colours create interesting effect in the decorations of faces and bodies.

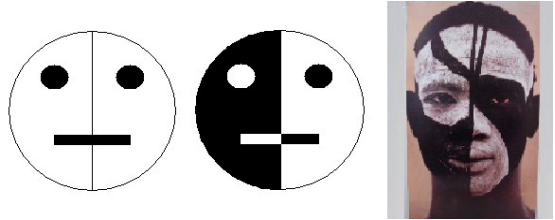

When seals are prepared to be used for stamp decorations on textiles, the creation of antisymmetric images with opposite colours is possible. In the case of those ancient seals, with images engraved in hard materials and used to stamp in clay, antisymmetric effects are very hard to create. Nevertheless we can ask ourselves what could be an alternative meaning of antisymmetry for such seals.

Let us consider the upper part of Figure 7, which shows a late bronze age seal from Aegeus with a lion and a horned quadruped [19]. As a referee told me, this seal could be viewed as representing a mirror antisymmetry, because the positions of animals are symmetric but animals are opposite to each other, in the sense that one is a predator and the other a prey. A proposal then for antisymmetry on seals is to have the mirror disposition of objects, but specific antithetic properties for them.

**Figure 7.** Seals of the late bronze age seal from Aegen area (adapted from [19,20]). The upper seal can be viewed as an antisymmetric seal.

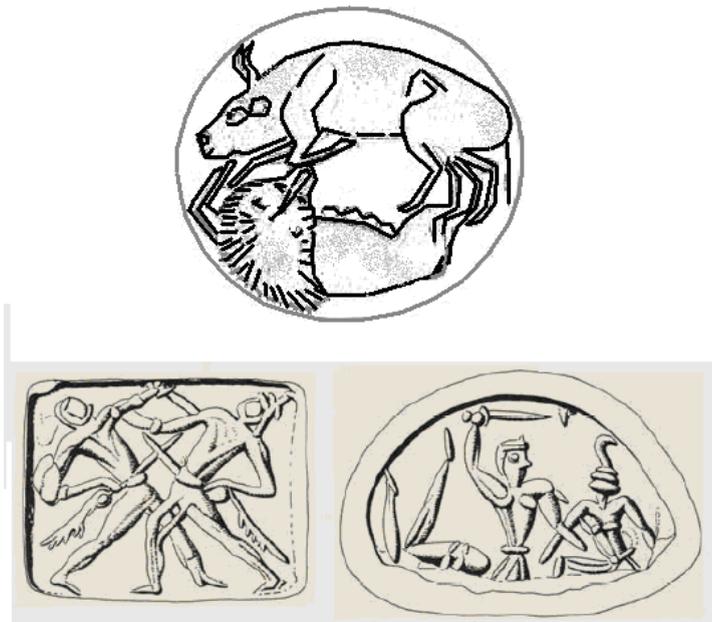

The hunting and fighting of wild animals and duels between heroes were a favourite themes of Mycenaean population [20], and then we find them often depicted on seals. Seals with a predator/prey couples were also quite common among Hittites [21]. The two seals in the lower part of Fig.7 are interesting for the arrangement of bodies in duels, bilateral on the left and possibly rotational on the right.

A discussion on Mycenaean seals is beyond the scope of this paper. Let us just note that during the Late Bronze Age, Mycenaean population adopted the sealing system of Minoan Crete. The most common shape of these seals was the lentoid, with a rich iconography. As the round shape of seals seems to favour images with rotational symmetries, in the case of lentoid seals, images often display the bilateral symmetry (Fig.9).

**Figure 9.** Mycenaean seals with bilateral symmetry.

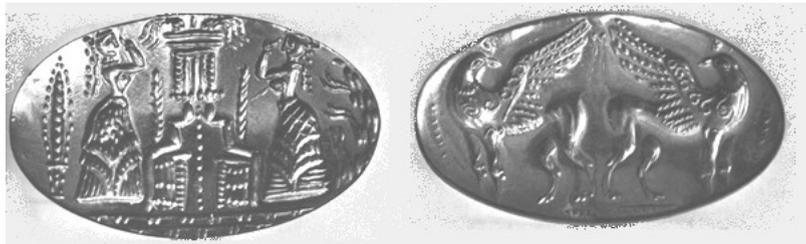

A symmetry can be exact, approximate, or broken. Exact means unconditionally valid, approximate means valid under certain conditions, but broken can mean different things, depending on the object considered and its context. The study of symmetry breaking in physics goes back to Pierre Curie. According to Curie, for the occurrence of a phenomenon in a certain medium, the original symmetry group of the medium must be lowered, that is broken, to the symmetry group of the phenomenon by the action of some cause. Then the symmetry breaking is what creates the phenomenon [22].

The breaking of a certain symmetry does not imply that no symmetry is present, but rather that the situation where this symmetry is broken is characterized by a lower symmetry than the situation where this symmetry is not broken.

Symmetry-breaking is a term used in the study of natural languages too: let me consider a phrase from a researcher in linguistic, who tells that "movement is a symmetry-breaking phenomenon" [23], and observe a beautiful Minoan stamp (Fig.8), with men fighting with bulls. This seal has four objects, and then in principle it could have a four-fold rotational symmetry: the fact that we have two bulls and two men reduces the symmetry to a two-fold one. Moreover, this symmetry is approximate. Note the impression of movement that this image is able to convey. In the symbolic language of seals, breaking the symmetry means an increase of movement, the creation of a phenomenon. In this case, it creates a scene from the bull-ring, with bulls careering in full gallop and athletes on them.

**Figure 8.** Men fighting with bulls, for Minoan area (adapted from [24]).

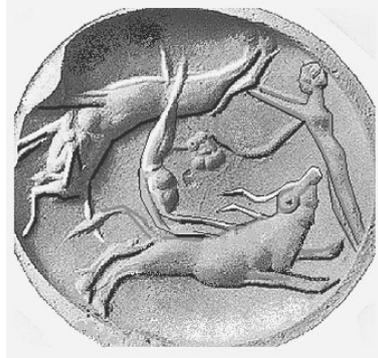

## 4. Egyptian scarab Seals

In the case of Egyptian seals we are in the lucky position that many catalogues of huge collections, for instance of British Museum or Egyptian Museum of Cairo [13-14], are in public domain. As a consequence, we can create a sort of statistics of images and check if symmetries are present and what are the preferred image arrangements. We find a relevant number of stamp seals with cords and coils, highly symmetric, mainly with bilateral symmetry, sometimes with two-fold rotational symmetry. Since scarabs have an elliptic basis, four-fold symmetry is not observed.

**Figure 10.** Images (a) with two-fold rotation symmetry with lions, scorpions and lizards, adapted from [14]. A seal (b) with lion and antelope, adapted from [13], and seals (c,d,e) of the Turin Egyptian Museum Collection, with bilateral and rotational symmetries.

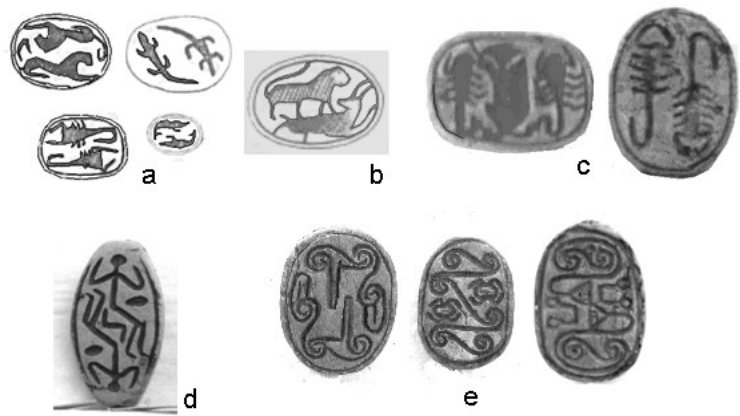

The two-fold symmetry seems to be not so frequent among Egyptian scarab seals depicting animals. It is more frequent when the image has coils and spirals. In catalogues, we find few scarabs with two symmetric animals (scorpions, lions, crocodiles or lizards, see Fig.10). It is not clear if these seals were unusual or not enough interesting to deserve a place in a museum collection, because the quality of image is sometimes poor. Since Scarabs were also amulet, we could guess for those seals with a rotational symmetry, an intent to increase the protection to the owner: the amulet protects in both cases, when it is upright or rotated upside down.

The collection of the Egyptian Museum of Torino possesses a unique seal with two men, may be twins (Fig.10 d). This is a cowroid seal, with lentoid shape, imitating the form of a cowry-shell, and can be dated to 2200-2040 BC, the first intermediate period of Egypt. The former scarabs with human figures were developed during this period, depicted in a linear style. Predator/prey or duels are subjects which are not interesting for Egyptian: the author found just a scarab in Ref.14 with a lion and an antelope. It is quite usual to find scarabs with a lion-hunting pharaoh, but they were amulets.

Three-fold, four-fold or greater order rotational symmetries are not observed on scarab seals, because it is difficult to arrange them in an elliptic shape. The collection of the Egyptian Museum of Turin has a seal with a five-fold rosette, but this is a button seals with a round shape. For what concerns the Egyptian knowledge of symmetry, B. Grünbaum claims that Egyptians appear to have missed the three-fold rotation [24]. May be this is not used in plane pattern decorations or for scarab seals, but the decoration of Egyptian faience pottery shows beautiful examples of three-fold symmetry (see Fig.11).

**Figure 11.** Egyptian faience with three-fold rotational symmetry, decorated with fishes and lotus flowers (Berlin, Ägyptisches Museum und Papyrussammlung).

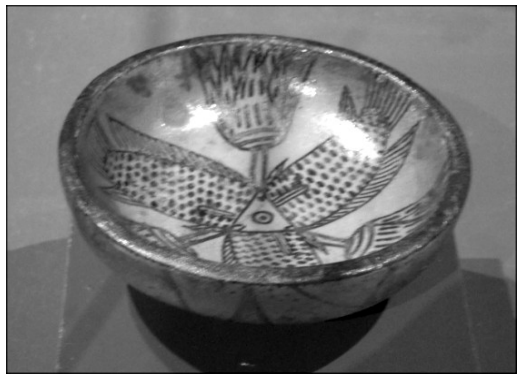

**5. Pre-Columbian seal symmetries**

There is an open problem concerning the pre-Colombian seals, "sellos" from Mexico or from other countries of Mesoamerica. This problem is if they were just devices to apply a design on textile and skin, or they were really seals, that is able to convey a divine authority upon the person stamped with them [10]. For sure, the use of seals to create coloured designs is part of the pre-Columbian tradition. But stamp seals were basically used for body decoration [10,26]: this idea is coming from the fact that there are ceramics with human shapes, showing bodies decorated with designs similar to those of the seals. If so, the use of seals for transient Tattooing confirm their symbolic values.

From scientific excavations, it is possible to conclude that stamps were part of the offerings in the burials as well as objects within housing areas. Pre-Columbian seals They have different sizes, with flat or cylindrical shape and designs which are simply figures or intricate patterns, as well as representations of the human figure and animal. The rectangular shape favours bilateral and two-fold rotational symmetries. Sometimes we see a sort of antisymmetry in those decorations with predator/prey couples (for instance the monkey and jaguar in Fig.12). Taino round stamp seals have four-fold rotational symmetry. Pre-Columbian seals are perhaps less known and studied, however they are of great importance in regard to cultural value.

**Figure 11.** Pre-Columbian seals with rectangular shape and bilateral (a) and rotational (b) symmetries. (c) shows a monkey and a jaguar (adapted from [26]). (e) and (f) are Taino round stamp seals.

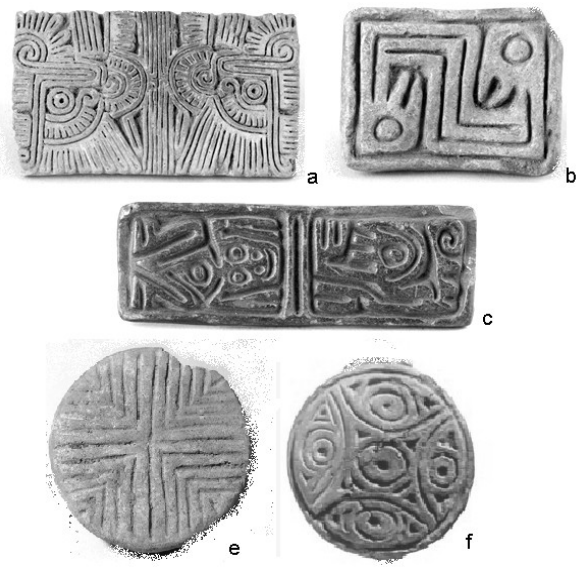

## 5. Conclusion

This paper discussed the presence of symmetries in images engraved on ancient seals. Bilateral and two-fold rotational symmetry are quite common, due to the fact that the shape of the base of stamp seals was often elliptic. Four-fold symmetries are sometimes observed in round stamp seals with geometric decorations.

The paper does not pretend to be exhaustive about all possibilities that can be observed on seals, because the research field and the related literature are huge. This paper would be a proposal for a new approach to seals analysis, in particular to stimulate the investigations on the meaning of antisymmetries and broken symmetries. Moreover, we have shown that there is a culture, the Indus Valley culture, that seems to avoid symmetries. Since it is common to find symmetries in human artifacts across time and cultures, a lack of symmetry is a very interesting phenomenon too.